\documentclass{amsart}%
\usepackage{graphicx}
\usepackage{amscd}
\usepackage{amsmath}
\usepackage{amsfonts}
\usepackage{amssymb}%
\providecommand{\U}[1]{\protect\rule{.1in}{.1in}}
\providecommand{\U}[1]{\protect\rule{.1in}{.1in}}
\providecommand{\U}[1]{\protect\rule{.1in}{.1in}}
\providecommand{\U}[1]{\protect\rule{.1in}{.1in}}
\newtheorem{theorem}{Theorem}
\theoremstyle{plain}

\newtheorem{corollary}{Corollary}
\newtheorem{definition}{Definition}

\newtheorem{lemma}{Lemma}

\newtheorem{remark}{Remark}
\numberwithin{equation}{section}

\begin{document}
\title{On the classification of quasi-homogeneous curves}
\author{Leonardo Meireles C\^{a}mara}
\address{Departamento de Matem\'{a}tica, CCE-UFES, Av. Fernando Ferrari 514, Campus de
Goiabeiras, Vit\'{o}ria-ES, 29075-910, Brazil.}
\email{camara@cce.ufes.br}
\thanks{This work was partially supported by CAPES-PROCAD grant n$^{\text{\d{o}}}$ 0007056}
\thanks{2000 Math. Subj. Class. Primary: 14H20. Secondary: 14B05, 14J15, 32S05.}

\begin{abstract}
We apply techniques of Holomorphic Foliations in the description of the
analytic invariants associated to germs of quasi-homogeneous curves in
$(\mathbb{C}^{2},0)$. As a consequence, we obtain an effective method to
determine whether two quasihomogeneous curves are analytically equivalent.

\end{abstract}
\maketitle

\section{Introduction}

The problem of the classification of germs of analytic plane curves has been
addressed by several authors since the XVII$^{\text{th}}$ century with
different methods (see for instance \cite{Be 06}, \cite{Br 66}, \cite{HefHer
2007}). In the present work, we study the problem of the analytic
classification of germs of singular curves with many branches from the
viewpoint of Holomorphic Foliations. This allows the use of geometrical
techniques including the blow-up and holonomy which are related to the study
of normal forms for quasi-homogeneous polynomials in two variables carried out
in \cite{Ca 2007}.

Next, we use the standard resolutions of theses singularities in order to
stratify them and thus identify the moduli space of each stratum. As a
consequence, our method provides an effective way to identify if two curves
are equivalent. Finally, we would like to remark that the analytic type of a
quasi-homogeneous curve is one of the invariants which determine the analytic
type of a foliation having such a curve as separatrix (cf. \cite{Ca 2007}).
Thus the present classification completes the classification of such germs of
complex analytic foliations.

\section{Preliminaries}

Let $C$ be a singular curve, $\pi:(\mathcal{M},D)\longrightarrow
(\mathbb{C}^{2},0)$ its standard resolution, i.e. the minimal resolution of
$C$ whose \textit{strict transform} $\widetilde{C}:=\pi^{-1}(C)-D$ \ is
transversal to the exceptional divisor $D=\pi^{-1}(0)$. A germ of holomorphic
function $f\in\mathbb{C}\{x,y\}$ is said to be \textit{quasi-homogeneous} if
there is a local system of coordinates in which $f$ can be represented by a
quasi-homogeneous polynomial, i.e. $f(x,y)=\sum_{ai+bj=d}a_{ij}x^{i}y^{j}$
where $a,b,d\in\mathbb{N}$. Let $M$ be a manifold and $M_{\Delta}%
(n):=\{(x_{1},\cdots,x_{n})\in M^{n}:x_{i}\neq x_{j}$ for all $i\neq j\}$. Let
$S_{n}$ denote the group of permutations of $n$ elements and consider its
action in $M_{\Delta}(n)$ given by $(\sigma,\lambda)\mapsto\sigma\cdot
\lambda=(\lambda_{\sigma(1)},\cdots,\lambda_{\sigma(n)})$. The quotient space
induced by this action is denoted by $Symm(M_{\Delta}(n))$. Now suppose a Lie
group $G$ acts in $M$ and let $G$ act in $M_{\Delta}(n)$ in the natural way
$(g,\lambda)=(g\cdot\lambda_{1},\cdots,g\cdot\lambda_{n})$ for every
$\lambda\in M_{\Delta}(n)$. Then the actions of $G$ and $S_{n}$ in $M_{\Delta
}(n)$ commute. Thus one obtains a natural action of $G$ in $Symm(M_{\Delta
}(n))$. Given $\lambda\in M_{\Delta}(n)$, denote its equivalence class in
$Symm(M_{\Delta}(n))/G$ by $[\lambda]$.

Let $C$ be a quasi-homogeneous curve determined by $f=0$, where $f$ is a
reduced polynomial. Then Lemma \ref{general decomp.} says that $f$ can be
(uniquely) written in the form%
\[
f(x,y)=x^{m}y^{k}%
{\displaystyle\prod\limits_{j=1}^{n}}
(y^{p}-\lambda_{j}x^{q})
\]
where $m,k\in\mathbb{Z}_{2}$, $p,q\in\mathbb{Z}_{+}$, $p\leq q$, $\gcd
(p,q)=1$, and $\lambda_{j}\in\mathbb{C}^{\ast}$ are pairwise distinct. In
particular $C$ has $n+m+k$ distinct branches. Since the exceptional divisor of
the standard resolution and the number of irreducible components are analytic
invariants of a germ of curve, then Lemmas \ref{companion fibration} and
\ref{companion 2} ensure that the triple $(p,q,n)$ is an analytic invariant of
the curve. Thus we have to consider the following three distinct cases:

\begin{enumerate}
\item[i)] $f(x,y)=x^{m}%
{\displaystyle\prod\limits_{j=1}^{n}}
(y-\lambda_{j}x)$ where $m\in\mathbb{Z}_{2}$, and $\lambda_{j}\in\mathbb{C}$.

\item[ii)] $f(x,y)=x^{m}%
{\displaystyle\prod\limits_{j=1}^{n}}
(y-\lambda_{j}x^{q})$ where $m\in\mathbb{Z}_{2}$, $q\in\mathbb{Z}_{+}$,
$q\geq2$ and $\lambda_{j}\in\mathbb{C}$.

\item[iii)] $f(x,y)=x^{m}y^{k}%
{\displaystyle\prod\limits_{j=1}^{n}}
(y^{p}-\lambda_{j}x^{q})$ where $m,k\in\mathbb{Z}_{2}$, $p,q\in\mathbb{Z}_{+}%
$, $2\leq p<q$, $\gcd(p,q)=1$, and $\lambda_{j}\in\mathbb{C}^{\ast}$.
\end{enumerate}

A quasi-homogeneous curve is said to be of \textit{type} $(1,1,n)$, $(1,q,n)$,
and $(p,q,n)$ respectively in cases i), ii), and iii).

\begin{theorem}
\label{main}The analytic moduli space of germs of quasi-homogeneous curves of
type $(p,q,n)$ are given respectively by\medskip

\begin{enumerate}
\item[i)] $\frac{Symm(\mathbb{P}_{\Delta}^{1}(n))}{PSL(2,\mathbb{C})}$ if
$(p,q)=(1,1)$;\bigskip

\item[ii)] $\mathbb{Z}_{2}\times\frac{Symm(\mathbb{C}_{\Delta}(n))}%
{Aff(\mathbb{C})}$ if $p=1$ and $q>1$;\bigskip

\item[iii)] $\mathbb{Z}_{2}\times\mathbb{Z}_{2}\times\frac{Symm(\mathbb{C}%
_{\Delta}^{\ast}(n))}{GL(1,\mathbb{C})}$ if $1<p<q$.
\end{enumerate}
\end{theorem}

\section{Quasi-homogeneous polynomials}

\subsection{Normal forms}

A\ quasi-homogeneous polynomial $f\in\mathbb{C}[x,y]$ is called
\textit{commode}\textbf{ }if its Newton polygon intersects both coordinate
axis. Further, notice that a polynomial in two variables $P\in\mathbb{C}[x,y]$
may be considered as a polynomial in the variable $y$ with coefficients in
$\mathbb{C}[x]$, i.e. $P\in(\mathbb{C}[x])[y]$. Let $\operatorname*{ord}_{y}P$
be the order of $P$ as a polynomial in $(\mathbb{C}[x])[y]$. Similarly let
$\operatorname*{ord}_{x}P$ be the order of $P$ as an element of $(\mathbb{C}%
[y])[x]$. Therefore, a quasi-homogeneous polynomial $P\in\mathbb{C}[x,y]$ is
commode if and only if $\operatorname*{ord}_{x}P=\operatorname*{ord}_{y}P=0$.
Next, we recall the general behavior of a quasi-homogeneous polynomial.

\begin{lemma}
\label{first decomp.}Let $P\in\mathbb{C}[x,y]$ be a quasi-homogeneous
polynomial, then it has a unique decomposition in the form
\[
P(x,y)=x^{m}y^{n}P_{0}(x,y)
\]
where $m,n\in\mathbb{N}$, $\lambda\in\mathbb{C}$, and $P_{0}$ is a commode
quasi-homogeneous polynomial.
\end{lemma}

\begin{proof}
Let $m:=\operatorname*{ord}_{x}P$ and $n:=\operatorname*{ord}_{y}P$. Clearly,
both $x^{m}$ and $y^{n}$ divide $P$. Hence $P$ can be written in the form
$P(x,y)=\sum_{ai+bj=d}a_{ij}x^{i}y^{j}$ where $i\geq m$ and $j\geq n$. Thus
$P(x,y)=x^{m}y^{n}P_{0}(x,y)$ where $P_{0}(x,y)=\sum_{ai^{\prime}+bj^{\prime
}=d^{\prime}}a_{i^{\prime}+m,j^{\prime}+n}x^{i^{\prime}}y^{j^{\prime}}$ and
$d^{\prime}:=d-am-bn$. Since $m=\operatorname*{ord}_{x}P$ and
$n=\operatorname*{ord}_{y}P$, then $\operatorname*{ord}_{x}P_{0}%
=0=\operatorname*{ord}_{y}P_{0}$. The result then follows directly from the
above remark.
\end{proof}

\begin{definition}
A commode polynomial $P\in\mathbb{C}[x,y]$ is called \textit{monic} in $y$ if
it is a monic polynomial in $(\mathbb{C}[x])[y]$.
\end{definition}

\begin{lemma}
\label{main decomp.}Let $P\in\mathbb{C}[x,y]$ be a commode quasi-homogeneous
polynomial, which is monic in $y$. Then $P$ can be written uniquely as
\[
P(x,y)=\prod_{\ell=1}^{k}(y^{p}-\lambda_{\ell}x^{q})\text{,}%
\]
where $\gcd(p,q)=1$ and $\lambda_{\ell}\in\mathbb{C}^{\ast}$.
\end{lemma}

\begin{proof}
First remark that any quasi-homogeneous polynomial can be written in the form
$P(x,y)=\sum_{pi+qj=m}a_{ij}x^{i}y^{j}$ where $p,q,m\in\mathbb{N}$ and
$\gcd(p,q)=1$. Since $P$ is commode, there are $i_{0},j_{0}\in\mathbb{N}$ such
that $qj_{0}=m$ and $pi_{0}=m$; in particular $k:=m/pq\in\mathbb{N}$.
Therefore $pi+qj=pqk$. Since $\gcd(p,q)=1$, then $q$ divides $i$ and $p$
divides $j$. If we let $i=qi^{\prime}$ and $j=pj^{\prime}$, then $pqi^{\prime
}+qpj^{\prime}=pqk$. Thus $P$ can be written in the form $P(x,y)=\sum
_{i+j=k}a_{qi,pj}x^{qi}y^{pj}$. Let $y=tx^{\frac{q}{p}}$, then the above
equation assumes the form $P(x,tx^{q/p})=x^{qk}\sum_{i+j=k}a_{qi,pj}t^{pj}$.
Now let $\{\lambda_{j}\}_{j=1}^{k}$ be the roots of the polynomial
$g(z)=\sum_{i+j=k}a_{qi,pj}z^{j}$, then
\begin{align*}
P(x,y)  &  =x^{qk}\prod_{\ell=1}^{k}(t^{p}-\lambda_{l})=x^{qk}\prod_{\ell
=1}^{k}(\frac{y^{p}}{x^{q}}-\lambda_{l})\\
&  =\prod_{\ell=1}^{k}(y^{p}-\lambda_{l}x^{q}).
\end{align*}

\end{proof}

\begin{lemma}
\label{general decomp.}Let $P\in\mathbb{C}[x,y]$ be a quasi-homogeneous
polynomial. Then $P$ can be written, uniquely, in the form%
\[
P(x,y)=\mu x^{m}y^{n}\prod_{\ell=1}^{k}(y^{p}-\lambda_{\ell}x^{q})
\]
where $m,n,p,q\in\mathbb{N}$, $\mu,\lambda_{\ell}\in\mathbb{C}^{\ast}$, and
$\gcd(p,q)=1$.
\end{lemma}

\begin{proof}
In view of Lemma \ref{first decomp.} and Lemma \ref{main decomp.}, it is
enough to remark that any commode quasi-homogeneous polynomial $P\in
\mathbb{C}[x,y]$ can be written uniquely as $P=\mu P_{0}$ where $P_{0}$ is
monic in $y$.
\end{proof}

\subsection{Resolution}

We recall the geometry of the exceptional divisor of the minimal resolution of
a germ of quasi-homogeneous curve.

A \textit{tree of projective lines} is an embedding of a connected and simply
connected chain of projective lines intersecting transversely in a complex
surface (two dimensional complex analytic manifold) with two projective lines
in each intersection. In fact, it consists of a pasting of Hopf bundles whose
zero sections are the projective lines themselves. A\textbf{\ }\textit{tree of
points} is any tree of projective lines in which a finite number of points is
discriminated. The above nomenclature has a natural motivation. In fact, it is
well know that one can assign to each projective line a point and to each
intersection an edge in other to form the \textit{weighted dual graph}. Two
trees of points are called \textit{isomorphic} if their weighted dual graph
are isomorphic (as graphs). \ It is well known that any germ of analytic curve
$C$ in $(\mathbb{C}^{2},0)$ has a standard resolution, which we denote by
$\widetilde{C}$. If the exceptional divisor of $\widetilde{C}$ has just one
projective line containing three or more singular points of $\widetilde{C}$,
then it is called the \textit{principal projective line} of $\widetilde{C}$
and denoted by $D_{\operatorname*{pr}(\widetilde{C})}$. A tree of projective
lines is called a \textit{linear chain} if each of its projective lines
intersects at most other two projective lines of the tree. A\ projective line
of a linear chain is called an \textit{end} if it intersects just another one
projective line of the chain.

\begin{lemma}
\label{companion fibration}Let $C$ be a commode quasi-homogeneous curve. Then
its standard resolution tree is a linear chain and its standard resolution
$\widetilde{C}$ intersects just one projective line of $D$, i.e. $C$ has one
of the following diagrams\ of resolution:%

{\includegraphics[
trim=0.000000in 0.000000in 0.004544in 0.000000in,
height=0.5725in,
width=2.3653in
]%
{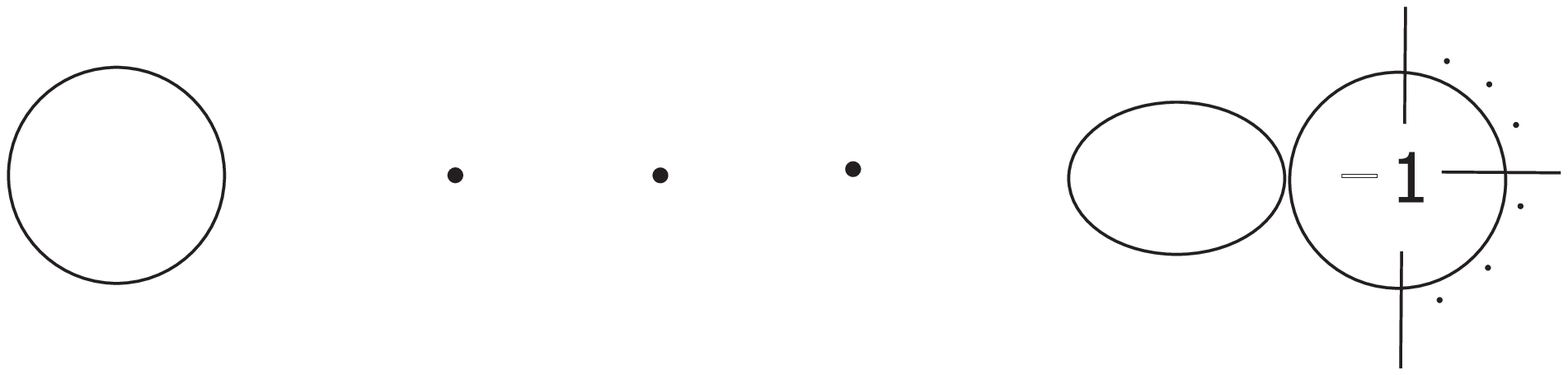}%
}
%

{\includegraphics[
trim=0.000000in 0.000000in -0.001454in 0.000000in,
height=0.5483in,
width=2.3367in
]%
{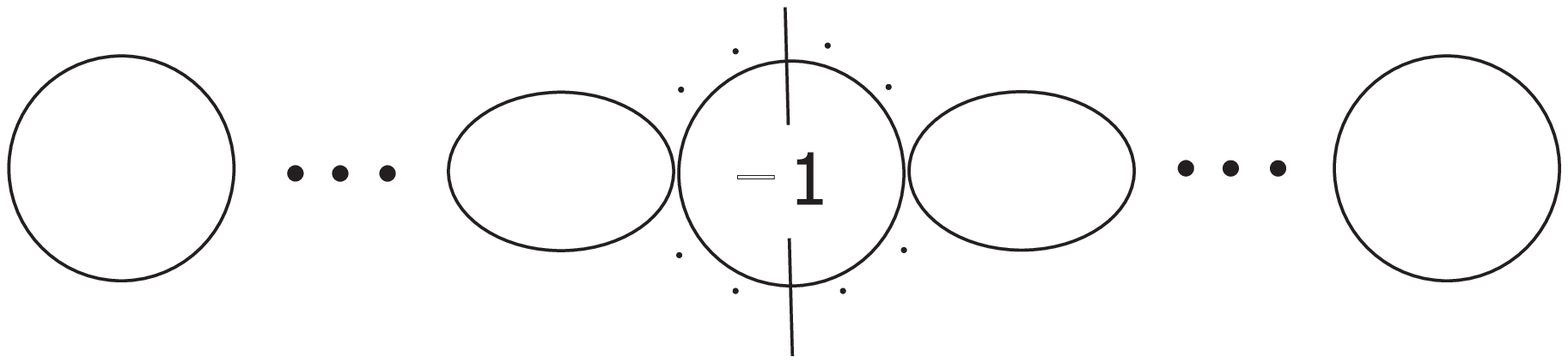}%
}

\end{lemma}

\begin{proof}
From Lemma \ref{main decomp.}, there is a local system of coordinates $(x,y)$
such that $C=f^{-1}(0)$ where $f(x,y)=\prod_{l=1}^{k}(y^{p}-\lambda_{j}x^{q})$
with $p<q$ and $\gcd(p,q)=1$. Since each irreducible curve $y^{p}-\lambda
_{l}x^{q}=0$ is a generic fiber of the fibration $\frac{y^{p}}{x^{q}}\equiv
const$, then it is resolved together with the fibration. After one blowup we
obtain:
\[%
\begin{array}
[c]{c}%
t^{p}/x^{q-p}\equiv const,\\
u^{q}y^{q-p}\equiv const.
\end{array}
\]
Since $p<q$, we have a singularity with holomorphic first integral at infinity
and a meromorphic first integral at the origin (as before). Going on with this
process, Euclid's algorithm assures that the resolution ends after the blowup
of a radial foliation. In particular, if $p=1$, then it is easy to see that
the\ principal projective line is transversal to just one projective line of
the divisor. Otherwise (i.e. if $p\neq1$) the singularity with meromorphic
first integral \textquotedblleft moves\textquotedblright\ to the
\textquotedblleft infinity\textquotedblright, i.e. it will appear in a corner
singularity. Then the principal projective line intersects exactly two
projective lines of the divisor.
\end{proof}

Let $\#irred(\widetilde{C})$ denote the number of irreducible components of
$\widetilde{C}$.

\begin{lemma}
\label{companion 2}Let $C$ be a non-commode quasi-homogeneous curve. Then its
minimal resolution tree is a linear chain having a principal projective line
such that $\#(\widetilde{\mathcal{C}}\cap D_{\operatorname*{pr}(\widetilde
{\mathcal{C}})})\leq$ $\#irred(\widetilde{\mathcal{C}})-1$. Further
$\widetilde{\mathcal{C}}\cap D_{j}=\emptyset$ whenever $D_{j}$ is neither the
principal projective line nor an end; i.e. $\mathcal{C}$ has one of the
following diagrams\ of resolution:%

{\includegraphics[
height=0.5621in,
width=2.3359in
]%
{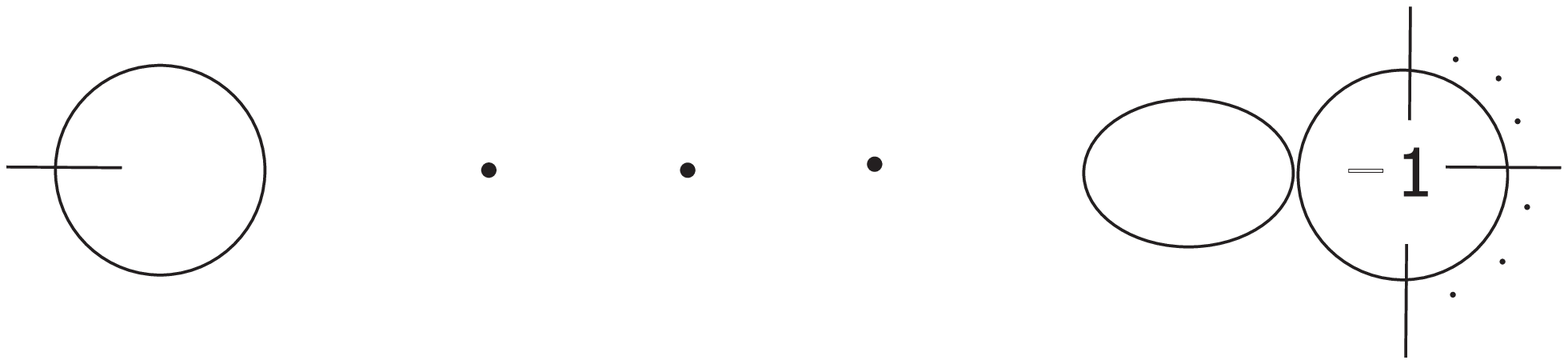}%
}
%

{\includegraphics[
trim=0.000000in 0.000000in -0.005302in 0.000000in,
height=0.5224in,
width=2.316in
]%
{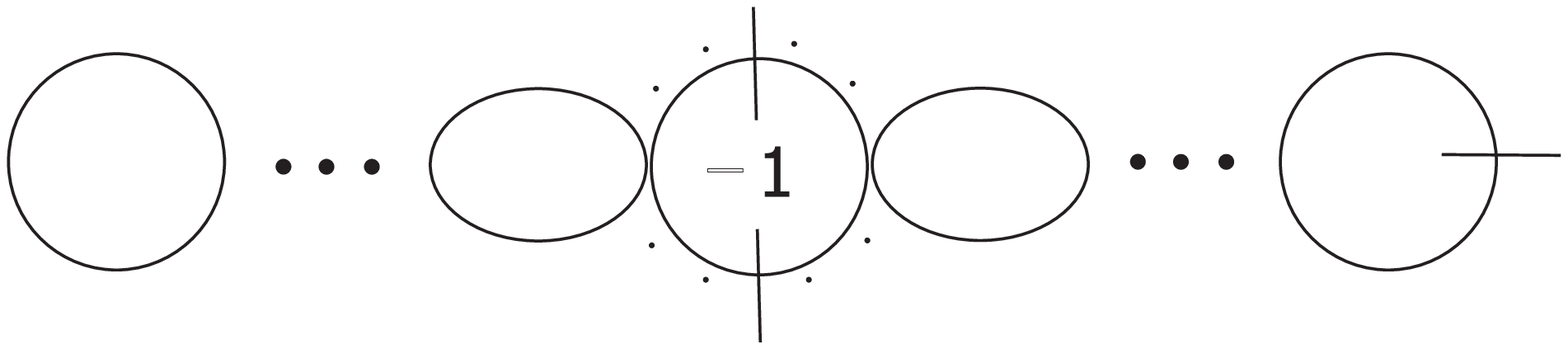}%
}
%

{\includegraphics[
trim=0.000000in 0.000000in 0.003124in 0.000000in,
height=0.512in,
width=2.3359in
]%
{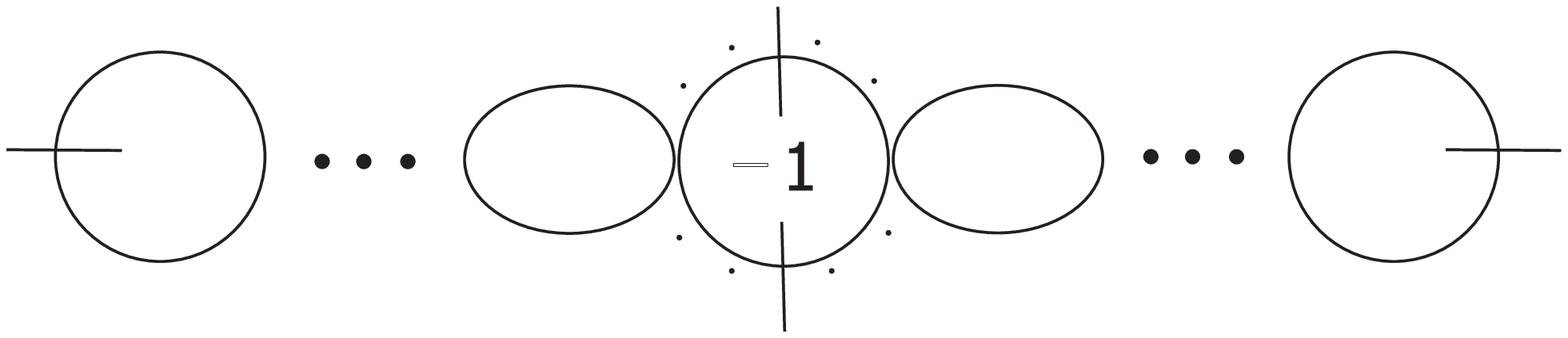}%
}

\end{lemma}

\begin{proof}
From Lemma \ref{general decomp.}, there is a local system of coordinates
$(x,y)$ such that $\mathcal{C}=f^{-1}(0)$ where $f(x,y)=\mu x^{m}y^{n}%
\prod_{l=1}^{k}(y^{p}-\lambda_{j}x^{q})$, $p<q$, and $\gcd(p,q)=1$. Since $\mu
x^{m}y^{n}$ is resolved after one blowup, then $f(x,y)$ is resolved together
with the fibration $\frac{y^{p}}{x^{q}}\equiv const$, as before. Then the
result follows from Lemma \ref{companion fibration}.
\end{proof}

\section{Quasi-homogeneous curves}

We consider each case separately and prove Theorem \ref{main} in a series of lemmas.

\subsection{Curves of type $(1,1,n)$.}

In this case the curve is given as the zero set of a polynomial of the form
$f(x,y)=x^{m}%
{\displaystyle\prod\limits_{j=1}^{n}}
(y-\lambda_{j}x)$ where $m\in\mathbb{Z}_{2}$, and $\lambda_{j}\in\mathbb{C}$;
in particular it is resolved after one blowup. Thus, given $\lambda
=(\lambda_{1},\cdots,\lambda_{n})\in\mathbb{P}_{\Delta}^{1}(n)$ we define
$f_{\lambda}(x,y)=x%
{\displaystyle\prod\limits_{j\neq i}}
(y-\lambda_{j}x)$ if $\lambda_{i}=\infty$ or $f_{\lambda}(x,y)=%
{\displaystyle\prod\limits_{j=1}^{n}}
(y-\lambda_{j}x)$ if $\lambda_{j}\neq\infty$ for all $j=1,\ldots,k$. We denote
the curve $f_{\lambda}=0$ by $C_{\lambda}$. Recall that the natural action of
$PSL(2,\mathbb{C})$ in $\mathbb{P}^{1}$ as the group of homographies induces a
natural action of $PSL(2,\mathbb{C})$ in $Symm(\mathbb{P}_{\Delta}^{1}(n))$.
Further, recall that the equivalence class of $\lambda\in$ $\mathbb{P}%
_{\Delta}^{1}(n)$ in $Symm(\mathbb{P}_{\Delta}^{1}(n))/PSL(2,\mathbb{C})$ is
denoted by $[\lambda]$.

\begin{lemma}
\label{(1,1)}Two homogeneous curves $C_{\lambda}$ and $C_{\mu}$ are
analytically equivalent if and only if $[\lambda]=[\mu]\in Symm(\mathbb{P}%
_{\Delta}^{1}(n))/PSL(2,\mathbb{C})$.
\end{lemma}

\begin{proof}
Suppose $C_{\lambda}$ and $C_{\mu}$ are analytically equivalent and let
$\Phi\in Dif(\mathbb{C}^{2},0)$ take $C_{\lambda}$ into $C_{\mu}$. Let
$\widetilde{\Phi}$ be the blowup of $\Phi$, then it takes the strict transform
of $C_{\lambda}$ into the strict transform of $C_{\mu}$. Blowing up
$f_{\lambda}$ and $f_{\mu}$ we obtain at once that the first tangent cones of
$C_{\lambda}$ and $C_{\mu}$ are respectively given by $\{\lambda_{1}%
,\cdots,\lambda_{n}\}$ and $\{\mu_{1},\cdots,\mu_{n}\}$. Therefore, there is
$\sigma\in S_{n}$ such that the M\"{o}bius transformation $\varphi=\left.
\widetilde{\Phi}\right\vert _{\mathbb{P}^{1}}$ satisfies $\mu_{\sigma
(j)}=\varphi(\lambda_{j})$ for all $j=1,\ldots,n$. In other words
$[\lambda]=[\mu]$. Conversely, suppose $[\lambda]=[\mu]$. Reordering the
indexes of $\{\mu_{1},\cdots,\mu_{n}\}$ we may suppose, without loss of
generality, that there is a M\"{o}bius transformation $\varphi(z)=\frac
{az+b}{cz+d}$, with $ad-bc=1$, such that $\mu_{j}=\varphi(\lambda_{j})$ for
all $j=1,\ldots,n$. Now consider the linear transformation
$T(x,y)=(dx+cy,bx+ay)$ with inverse $T^{-1}(x,y)=(ax-cy,-bx+dy)$. Then a
straightforward calculation shows that $f_{\lambda}=\alpha\cdot T^{\ast}%
f_{\mu}$ where $\alpha\in\mathbb{C}^{\ast}$. Thus $C_{\lambda}$ is
analytically equivalent to $C_{\mu}$, as desired.
\end{proof}

\begin{remark}
\label{three branches}Recall that for any three distinct points $\{\lambda
_{1},\lambda_{2},\lambda_{3}\}\subset\mathbb{P}^{1}$ there is a M\"{o}bius
transformation $\varphi$ such that $\varphi(0)=\lambda_{1}$, $\varphi
(1)=\lambda_{2}$ and $\varphi(\infty)=\lambda_{3}$.
\end{remark}

As a straightforward consequence of Lemma \ref{(1,1)} and Remark
\ref{three branches} one has:

\begin{corollary}
Let $\lambda,\mu\in\mathbb{P}_{\Delta}^{1}(n)$ with $n\leq3$. Then
$C_{\lambda}$ and $C_{\mu}$ are analytically equivalent.
\end{corollary}

\subsection{Curves of type $(1,q,n)$, $q\geq2$.}

In this case, the curve is given as the zero set of a polynomial of the form
$f_{m,\lambda}(x,y)=x^{m}%
{\displaystyle\prod\limits_{j=1}^{n}}
(y-\lambda_{j}x^{q})$ where $m\in\mathbb{Z}_{2}$, $q\in\mathbb{Z}_{+}$,
$q\geq2$, and $\lambda_{j}\in\mathbb{C}$. Thus given $m\in\mathbb{Z}_{2}$ and
$\lambda=(\lambda_{1},\cdots,\lambda_{n})\in\mathbb{C}_{\Delta}(n)$, we denote
a curve of type $(1,q,n)$ by $C_{m,\lambda}$ if it is given as the zero set of
$f_{m,\lambda}$. Recall that the group of affine transformations of
$\mathbb{C}$, denoted by $Aff(\mathbb{C})$, acts in a natural way in
$Symm(\mathbb{C}_{\Delta}(n))$. Further, recall that the equivalence class of
$\lambda\in$ $\mathbb{C}_{\Delta}(n)$ in $Symm(\mathbb{C}_{\Delta
}(n))/Aff(\mathbb{C})$ is denoted by $[\lambda]$.

\begin{lemma}
\label{(1,q)}Two homogeneous curves $C_{m,\lambda}$ and $C_{m,\mu}$ are
analytically equivalent if and only if $[\lambda]=[\mu]\in Symm(\mathbb{C}%
_{\Delta}(n))/Aff(\mathbb{C})$.
\end{lemma}

\begin{proof}
Suppose $\Phi\in Diff(\mathbb{C}^{2},0)$ is an equivalence between
$C_{m,\lambda}$ and $C_{m,\mu}$. From the proof of Lemma
\ref{companion fibration}, both curves are resolved after $q$ blowups.
Further, after $q-1$ blowups $\Phi$ will be lifted to a local conjugacy
$\Phi^{(q-1)}$ between the germs of curves given in local coordinates $(x,y)$
respectively by $p_{\lambda}(x,y)=x%
{\displaystyle\prod\limits_{j=1}^{n}}
(y-\lambda_{j}x)$ and $p_{\mu}(x,y)=x%
{\displaystyle\prod\limits_{j=1}^{n}}
(y-\mu_{j}x)$ where $(x=0)$ is the local equation of the exceptional divisor
$D^{(q-1)}$. Let $\pi$ denote a further blowup given in local coordinates by
$\pi(t,x)=(x,tx)$ and $\pi(u,y)=(u,uy)$, and $\Phi^{(q)}$ be the map obtained
by the lifting of $\Phi^{(q-1)}$ by $\pi$. Further, let $\varphi=\left.
\Phi^{(q)}\right\vert _{D_{q}}$ where $D_{q}=\pi^{-1}(0)$. Since $\Phi^{(q)}$
preserves the irreducible components of $\pi^{\ast}(D^{(q-1)})$, then
$\varphi(t)=\Phi^{(q)}(t,0)$ is a homography fixing $\infty$ and conjugating
the first tangent cones of $p_{\lambda}=0$ and $p_{\mu}=0$ respectively. Thus
$[\lambda]=[\mu]\in Symm(\mathbb{C}_{\Delta}(n))/Aff(\mathbb{C})$. Conversely,
(reordering the indexes of $\mu$, if necessary) suppose there is
$\varphi(z)=az+b\in Aff(\mathbb{C})$ such that $\mu_{j}=\varphi(\lambda_{j})$
for all $j=1,\ldots,n$, and let $T(x,y)=(x,ay+bx^{q})$. Then a straightforward
calculation shows that $f_{m,\lambda}=\alpha\cdot T^{\ast}f_{m,\mu}$ where
$\alpha\in\mathbb{C}^{\ast}$. Thus $C_{m,\lambda}$ and $C_{m,\mu}$ are
analytically equivalent, as desired.
\end{proof}

As a straightforward consequence of Lemma \ref{(1,q)} and Remark
\ref{three branches} one has:

\begin{corollary}
Let $\lambda,\mu\in\mathbb{C}_{\Delta}(n)$ with $n\leq2$. Then $C_{m,\lambda}$
and $C_{m,\mu}$ are analytically equivalent.
\end{corollary}

\subsection{Curves of type $(p,q,n)$, $2\leq p<q$.}

In this case, the curve is given as the zero set of a polynomial of the form
$f_{m,k,\lambda}(x,y)=x^{m}y^{k}%
{\displaystyle\prod\limits_{j=1}^{n}}
(y^{p}-\lambda_{j}x^{q})$ where $m,k=0,1$, $p,q\in\mathbb{Z}_{+}$, $2\leq
p<q$, and $\lambda_{j}\in\mathbb{C}^{\ast}$. Thus given $\lambda=(\lambda
_{1},\cdots,\lambda_{n})\in\mathbb{C}_{\Delta}^{\ast}(n)$ we denote a curve of
type $(p,q,n)$ by $C_{m,k,\lambda}$ if it is given as the zero set of
$f_{m,k,\lambda}(x,y)$. Recall that the group of linear transformations of
$\mathbb{C}$, denoted by $GL(1,\mathbb{C})$, acts in a natural way in
$Symm(\mathbb{C}_{\Delta}^{\ast}(n))$. Further, recall that the equivalence
class of $\lambda\in$ $\mathbb{C}_{\Delta}^{\ast}(n)$ in $Symm(\mathbb{C}%
_{\Delta}^{\ast}(n))/GL(1,\mathbb{C})$ is denoted by $[\lambda]$.

\begin{lemma}
\label{(p,q)}Two homogeneous curves $C_{m,k,\lambda}$ and $C_{m,k,\mu}$ are
analytically equivalent if and only if $[\lambda]=[\mu]\in Symm(\mathbb{C}%
_{\Delta}^{\ast}(n))/GL(1,\mathbb{C})$.
\end{lemma}

\begin{proof}
First recall from the proof of Lemma \ref{companion fibration} that
$C_{m,k,\lambda}$ is resolved after $N$ blowups, where $N$ depends on the
Euclid's division algorithm between $q$ and $p$. Further, in the $(N-1)^{th}$
step we have to blowup a singularity given in local coordinates $(x,y)$ as the
zero set of the polynomial $g_{\lambda}(x,y)=xy%
{\displaystyle\prod\limits_{j=1}^{n}}
(y-\lambda_{j}x)$. Therefore, if $\Phi\in Diff(\mathbb{C}^{2},0)$ is an
equivalence between $C_{m,k,\lambda}$ and $C_{m,k,\mu}$ and $\Phi^{(N-1)}$ is
its lifting to the $(N-1)^{th}$ step of the resolution, then it conjugates the
germs of curves given in local coordinates $(x,y)$ respectively by
$p_{\lambda}(x,y)=xy%
{\displaystyle\prod\limits_{j=1}^{n}}
(y-\lambda_{j}x)$ and $p_{\mu}(x,y)=xy%
{\displaystyle\prod\limits_{j=1}^{n}}
(y-\mu_{j}x)$ where $(x=0)$ and $(y=0)$ are local equations for the
exceptional divisor $D^{(N-1)}$. Let $\pi$ denote the final blowup of the
resolution given in local coordinates by $\pi(t,x)=(x,tx)$ and $\pi
(u,y)=(u,uy)$, and $\Phi^{(N)}$ be the map obtained by the lifting of
$\Phi^{(N-1)}$ by $\pi$. Further let $\varphi=\left.  \Phi^{(N)}\right\vert
_{D_{N}}$ where $D_{N}=\pi^{-1}(0)$. Since $\Phi^{(N)}$ preserves the
irreducible components of $\pi^{\ast}(D^{(q-1)})$, then $\varphi(t)=\Phi
^{(q)}(t,0)$ is a homography fixing $0$ and $\infty$, and conjugating the
first tangent cones of $p_{\lambda}=0$ and $p_{\mu}=0$ respectively. Thus
$[\lambda]=[\mu]\in Symm(\mathbb{C}_{\Delta}^{\ast}(n))/GL(1,\mathbb{C})$.
Conversely, (reordering the indexes of $\mu$, if necessary) suppose there is
$\varphi(z)=az\in GL(1,\mathbb{C})$ such that $\mu_{j}=\varphi(\lambda_{j})$
for all $j=1,\ldots,n$, and let $T(x,y)=(x,\sqrt[p]{a}y)$. Then a
straightforward calculation shows that $f_{m,\lambda}=\alpha\cdot T^{\ast
}f_{m,\mu}$ where $\alpha\in\mathbb{C}^{\ast}$. Thus $C_{m,\lambda}$ and
$C_{m,\mu}$ are analytically equivalent, as desired.
\end{proof}

As a straightforward consequence of Lemma \ref{(p,q)} and Remark
\ref{three branches} one has:

\begin{corollary}
Let $\lambda,\mu\in\mathbb{C}_{\Delta}^{\ast}(1)$, then $C_{m,k,\lambda}$ and
$C_{m,k,\mu}$ are analytically equivalent.
\end{corollary}

\section{Resolution and factorization}

We study the relationship between the resolution tree and the factorization of
a quasi-homogeneous polynomial. We use the resolution in order to study the
equivalence between two quasi-homogeneous polynomials.

First recall that a quasi-homogeneous polynomial split uniquely in the form
$P=x^{m}y^{n}P_{0}$ where $P_{0}$ is a commode quasi-homogeneous polynomial.
In particular $P$ and $P_{0}$ share the same resolution process.

\begin{corollary}
Let $P\in\mathbb{C}[x,y]$ be a commode quasi-homogeneous polynomial with the
weights $(p,q)$, where $\gcd(p,q)=1$. Let $q_{j}=s_{j}p_{j}+r_{j}$,
$j=1,\ldots,m$, be the Euclid's algorithm of $(p,q)$, where $q_{1}:=q$,
$p_{1}:=p$, $q_{j+1}:=p_{j}$, and $p_{j+1}:=r_{j}$ for all $j=1,\ldots,m-1$.
Then the exceptional divisor of its minimal resolution is given by a linear
chain of projective lines, namely $D=\cup_{j=1}^{n}D_{j}$, whose
self-intersection numbers are given as follows:

\begin{enumerate}
\item If $m=1$, then
\[
c_{1}(D_{j})=\left\{
\begin{array}
[c]{ll}%
-1 & \text{if }j=s_{1}\text{;}\\
-2 & \text{otherwise}%
\end{array}
\right.
\]

\item If $m=2\alpha+1$, $\alpha\geq1$, then
\[
c_{1}(D_{j})=\left\{
\begin{array}
[c]{ll}%
-(s_{2k}+2) & \text{if }j=s_{1}+s_{3}+\cdots+s_{2k-1}\text{, }k=1,\ldots
,\alpha\text{;}\\
-1 & \text{if }j=s_{1}+s_{3}+\cdots+s_{2\alpha+1}\text{;}\\
-(s_{2k+1}+2) & \text{if }j=m-(s_{2}+s_{4}+\cdots+s_{2k})+1\text{, }%
k=1,\ldots,\alpha-1\text{;}\\
-(s_{2\alpha+1}+1) & \text{if }j=m-(s_{2}+s_{4}+\cdots+s_{2\alpha}%
)+1\text{;}\\
-2 & \text{otherwise}%
\end{array}
\right.
\]

\item If $m=2\alpha$,  $\alpha\geq1$, then
\[
c_{1}(D_{j})=\left\{
\begin{array}
[c]{cl}%
-(s_{2k}+2) & \text{if }j=s_{1}+s_{3}+\cdots+s_{2k-1}\text{, }k=1,\ldots
,\alpha-1\text{;}\\
-(s_{2\alpha}+1) & \text{if }j=s_{1}+s_{3}+\cdots+s_{2\alpha-1}\text{;}\\
-(s_{2k+1}+2) & \text{if }j=m-(s_{2}+s_{4}+\cdots+s_{2k})+1\text{, }%
k=1,\ldots,\alpha-1\text{;}\\
-1 & \text{if }j=m-(s_{2}+s_{4}+\cdots+s_{2\alpha})+1\text{;}\\
-2 & \text{otherwise.}%
\end{array}
\right.
\]
Finally, if $C$ is given by $f=0$ where $f(x,y)=x^{m}y^{n}%
{\displaystyle\prod\limits_{j=1}^{k}}
(y^{p}-\lambda_{j}x^{q})$, then a representative of $[\lambda]$ is determined
by the intersection of the strict transform of $C$\ with the exceptional
divisor $D$.
\end{enumerate}
\end{corollary}

\begin{proof}
The proof shall be performed by induction on $m$, the length of the Euclidean
algorithm. In order to better understand the arguments, the reader have to
keep in mind the proof of Lemma \ref{companion fibration}. From Lemma
\ref{main decomp.}, we may suppose without loss of generality that $P$ can be
written in the form $P(x,y)=%
{\displaystyle\prod\limits_{j=1}^{k}}
(y^{p}-\lambda_{j}x^{q})$. First remark that if $m=1$ then $p=1$. Thus we
prove the statement for $m=1$ by induction on $q$. For $q=1$ the result is
easily verified after one blowup. Now suppose the result is true for all
$q\leq q_{0}-1$. Then after one blowup $\pi(t,x)=(x,tx)$, $\pi(u,y)=(uy,y)$,
$P$ is transformed into $\pi^{\ast}P(t,x)=x%
{\displaystyle\prod\limits_{j=1}^{k}}
(t-\lambda_{j}x^{q-1})$. Thus the result follows for $m=1$ by induction on
$q$. Suppose the result is true for all polynomials whose pair of weights have
Euclid's algorithm length less than $m$, and let $(p,q)$ has length $m$. Since
$p_{j}=s_{j}q_{j}+r_{j}$, $j=1,\ldots,m$, is the Euclid's algorithm of
$(p,q)$, then $p_{j}=s_{j}q_{j}+r_{j}$, $j=2,\ldots,m$, is the Euclid's
algorithm of $(p_{2},q_{2})$. In particular the Euclid's algorithm of
$(p_{2},q_{2})$ has length $m-1$. Reasoning in a similar way as in the case
$m=1$, we have after $s_{1}$ blowups a linear chain of projective lines
$\cup_{j=1}^{s_{1}}D_{j}^{(1)}$ such that $c_{1}(D_{j}^{(1)})=-2$ for all
$j=1,\ldots,s_{1}-1$ and $c_{1}(D_{s_{1}}^{(1)})=-1$. Further, the strict
transform of $P=0$ is given by the zero set of the polynomial $\widetilde
{P}(t,x)=%
{\displaystyle\prod\limits_{j=1}^{k}}
(t^{p_{1}}-\lambda_{j}x^{r_{1}})=$ $\lambda_{1}\cdots\lambda_{k}%
{\displaystyle\prod\limits_{j=1}^{k}}
(x^{p_{2}}-\lambda_{j}t^{q_{2}})$ where the local equation for $D_{s_{1}%
}^{(1)}$ is $(x=0)$. The first statement thus follows by the induction
hypothesis. The last statement comes immediately from the above reasoning.
\ For the above induction arguments ensure that the strict transform of $P$
assume the form $\widetilde{P}=0$, with $\widetilde{P}(x,y)=%
{\displaystyle\prod\limits_{j=1}^{k}}
(y-\lambda_{j}x)$, just before the last blowup.
\end{proof}

The above Corollary gives an effective way to compute the relatively prime
weights of a quasi-homogeneous polynomials and further an easy way to
determine a normal form for a quasi-homogeneous function from the dual
weighted tree of its standard resolution. In particular it shows how to split
a quasi-homogeneous polynomial into irreducible components from the dual
weighted tree of its minimal resolution.


\begin{thebibliography}{99}                                                                                               %


\bibitem {Be 06}Berzolari, L. \textit{Allgemeine Theorie der h\"{o}hernen
alagebrischen Kurven.} Enzyklop\"{a}die der math. Wissenschaften, Bd. III 2,
1, 313-455 (1906);

\bibitem {Br 66}Brieskorn, E. \textit{Beispiele zur Differentietopologie von
Singularit\"{a}ten}. Invent. Math. \textbf{2} (1966), 1-14;

\bibitem {Ca 2007}C\^{a}mara, L. M. \textit{On the moduli space of
quasihomogeneous foliations on} $(\mathbb{C}^{2},0)$\textit{. }IMPA preprint
serie D 42/2007. http://www.preprint.impa.br/Shadows/SERIE\_D/2007/42.html.

\bibitem {De 78}Delorme, C. \textit{Sur les Modules des Singularit\'{e}s de
Courbes Planes}. Bull. Soc. Math. France 106, 417-446 (1978);

\bibitem {Ebey 65}Ebey, S. \textit{The classification of singular points of
algebraic curves}. Trans. Amer. Math. Soc. \textbf{118}, 454-471 (1965).

\bibitem {Greuel 97}Greuel, G.-M., Hertling, C. \& Pfister, G. \textit{Moduli
spaces of semiquahomogeneous singularities with fixed principal part}. J.
Algebraic Geom. \textbf{6} (1997), 169-199;

\bibitem {HefHer 2007}Hefez, A. \& Hernandes, M. E. \textit{The analytic
classification of plane branches}. arXiv:0707.4502v1 [math.AG] 30 Jul 2007.

\bibitem {LinZai 83}Lin, V. I. \& Zaidenberg, M. \textit{An irreducible simply
connected algebraic curve in }$\mathbf{C}^{2}$\textit{ is equivalent to a
quasihomogeneous curve}. Dokl. Akad. Nauk. SSSR 271(1983), no5, 1048-1052.

\bibitem {Mil 68}Milnor, J.\textit{\ Singular points of complex
hypersurfaces.} Ann. Math. Studies \textbf{62}, Princeton University Press,
Princeton, N.Y., 1968.

\bibitem {MilnOlr 70}Milnor, J. \& Orlik, P.\textit{\ Isolated singularities
defined by homogeneous polynomials.} Topology \textbf{9} (1970), 385-393.

\bibitem {Po 74}Poincar\'{e}, H.\textit{\ Note sur les propri\'{e}t\'{e}s des
fonctions d\'{e}finies par des \'{e}quations diff\'{e}rentielles}, J. Ec.
Pol., \textbf{45}$^{\text{e}}$ cahier (1878), 13-26.

\bibitem {Re 54}Reeve, J. E. \textit{A\ summary of results in the topological
classification of plane algebroid singularities}. Universit\`{a} e Politecnico
di Torino, Rendiconti del Seminaro Matematico 1954/1955.

\bibitem {Sai 71}Saito, K. \textit{Quasihomogen isoliere Singularit\"{a}ten
von Hyperfl\"{a}chen}. Invent. math. 14(1971), 123-142.

\bibitem {Se 68}Seindenberg, A.\textit{\ Reduction of the singularities of the
differentiable equation }$Ady=Bdx$, Amer. J. Math \textbf{90} (1968), 248-269.

\bibitem {Za 66}Zariski, O. \textit{On the topology of algebroid
singularities,} Amer. J. Math. \textbf{54} (1932), 453-465.

\bibitem {Za 72}Zariski, O. \textit{Characterization of plane algebroid curves
whose module of differentials has maximum torsion}. Proceedings of the
National Academy of Science, \textbf{56 (3)} (1966), 781-786.

\bibitem {Za 73}Zariski, O. \textit{Le Probl\'{e}me des Modules pour les
Branches Planes}. Cours donn\'{e} au Centre de Math\'{e}matiques de
L'\'{E}cole Polytechnique. Nouvelle \'{e}d. revue par l'auteur. R\'{e}dig\'{e}
par Fran\c{c}ois Kimety et Michel Merle. Avec un appendice de Bernard
Teissier. Paris, Hermann (1986). English translation by Ben Lichtin: The
Moduli Problem for Plane Branches. University Lecture Series, Vol. 39, AMS (2006).
\end{thebibliography}
\end{document}